\documentclass{article}
\usepackage{amsfonts}
\newtheorem{theorem}{Theorem}[section]
\newtheorem{pro}{Proposition}[section]
\newtheorem{lemma}{Lemma}[section]
\newtheorem{sublemma}{Sublemma}[section]

\newcommand{\proof}{\noindent{\bf Proof}:\hspace*{3mm}}
\newcommand{\QED}{\begin{flushright} QED\end{flushright}}
\newcommand{\udn}{{\tilde u}_{d', N'}}
\newcommand{\da}{\frac{\partial}{\partial a}}
\newcommand{\dth}{\frac{\partial}{\partial \theta}}
\newcommand{\dyi}{\frac{\partial}{\partial y_i}}
\newcommand{\ds}{\frac{\partial}{\partial s}}
\newcommand{\dt}{\frac{\partial}{\partial t}}
\newcommand{\bp}{\bar\partial}
\newcommand{\sz}{S^2\setminus\{z_-, z_+\}} 

\title{  Fredholm theory  of the linearized ${{\bar \partial }}$-operator
 and additivity of its index }
\author{Gang Liu}
\begin{document}
\maketitle

\section{Introduction}

 In [L1], we established 
 the compactification of the moduli space  
of the  pseudo-holomorphic maps connecting the two ends of  the 
symplectization  of a compact  contact manifold. Several new  phenomena
concerning  the behavior of the boundary of 
 the compactification were described there.  

The first of them is that each time when a family of 
connecting pseudo-holomorphic maps bubbles off a bubble,
there is a new component of the target lying on the "left'
the original one such that the image of the bubble lies on
this new component of the target. Moreover, each time when
a new top bubble appears,  
there is also a new principal component lying on
the new component of the target. Therefore, even in the simplest case when
  such a 
family of connecting maps develops only one bubble, the domain of
the limit map necessarily has three components. This is quite different
from  what happens in Gromov-Floer theory. Note that in the bubbling
process above, both the domain and the target of the connecting map
split (degenerate) into different components. Furthermore, the rates of 
the two types of degeneration are independent of each other.
Closely related to this is that the 
${\bf R}$-symmetry of the target splits into a two-dimensional or
multi-dimensional symmetries during the bubbling or splitting process.

The second of new phenomena we found is that in order to get
a correct definition of a stable map and use it to get the desired
compactification, it is necessary to count  the ${\bf R}^1$-symmetry
of a component of the target as many times as the number of the connected
components of the stable map in the component of the target. 

These new phenomena, once carefully analyzed, lead to a rather
different picture on  the behavior of the  boundary
of the compactification of the moduli space of the
connecting pseudo-holomorphic maps from  what we used to believed.
It opens a door to various new constructions in contact geometry,
some of which were  outlined in [L1]. 

In fact, the results of [L1] lead to the following alternatives
depending on whether or not the additivity of the index formula
on the linearized ${\bar {\partial}}$-operator is applicable
in the contact case. To motivate the work of this paper,
we now  explain this alternative.   

Recall that since in the symplectization the symplectic form is exact,
 each top  bubble in the bubbling 
necessarily has non-removable singularity at infinity, and along the end at 
infinity,
the bubble is convergent to some closed orbit of the Reeb field 
(contact field)
of the
 contact 
manifold. Similar to the usual Floer homology, a family of connecting
pseudo-holomorphic maps between two closed orbits at the two ends may
split into a family of broken connecting maps. 
Same as [L1], we are going to assume throughout this paper
that the contact structure and the contact $1$-form are generic in the
sense that along the normal direction of each closed orbit of the
contact field, the linearized Poincare returning map does not have any
eigen vector of eigen value $1$.  This implies that the set of 
unparameterized closed  orbits is discrete.
 As it is well known that
the parameterized closed orbits of the contact field can be thought as
the critical points for the action functional
$a_{\lambda}: {\cal L}\rightarrow {\bf R}$ given by
the formula 
$$a_{\lambda}(z)= \int (z^* \lambda)(t) dt$$ for any $z\in {\cal L(M)},$ 
where ${\cal L(M)}$ is the free loop space of parameterized loops in the
contact manifold $M$ with the contact form $\lambda$.
Our assumption implies that the set of critical points is a discrete
union
of one-dimensional critical manifolds, each of them is a circle of
length $c=\int  x^*{\lambda} dt,$ where $x(t)$ is a closed orbit in
some critical manifold. This brings us to a situation of a
particular  Bott-type
Morse theory with the infinite dimensional action functional
$a_{\lambda}.$  The index homology in contact geometry introduced
in [L3] is a certain kind of Morse-Floer homology for the action
functional.

To define the index homology, we first define its chain complex as 
follows. Let $\{x\}$ be the equivalent class of the parameterized
 closed orbit
$x$, where any two of parameterized closed orbits are said to be equivalent
if they differ by a linear parameterization.  The chain complex is defined 
to be

$$ IC= IC (\lambda)= \oplus_{x\in {\cal P}} {\bf Q} \{x\},$$
where ${\cal P}$ is the set of parameterized closed  orbits of the contact
field.

 We define the boundary operator $D$ of the chain 
complex as follows.
Given $\{x_{+}\}\in IC$,  
$ D(\{x_{+}\}) = \Sigma_{x_{-}\in {\cal P}} <\{x_{-}\}, \{x_{+}\}>\{x_{-}\}$,
where $<\{x_{-}\}, \{x_{+}\}>=
\# \{ \cup_{x_{-}\in {\cal P}} {\cal M}( x_{-}, x_{+}, J)\}.$
 Here  ${\cal M}( x_{-}, x_{+}, J) $ is the moduli space of all
equivalent classes of 
${\tilde J}$-holomorphic maps in the symplectization
${\tilde M}= M\times {\bf R}$ connecting the two closed orbits $x'_{-}\in
\{x_{-}\}$ and  $x'_{+}\in
\{x_{+}\}$  in the two ends $M_{-}$ and $M_{+}$ of ${\tilde M}$. In
other words, we define the boundary operator by  counting discrete
connecting maps. Of course, to get a precise definition one needs to
use the corresponding virtual moduli cycles to replace the moduli space
used here ( for the construction of the virtual moduli cycles in the contact
case, see [L2]). 

To prove the fact that   $D^2=0$,  it is crucial to show
that the the virtual co-dimension of the moduli space 
${\cup}_{x\in {\cal P},\{x\}\not =\{x_{-}\}, \{x\}\not =\{x_{+}\}}
{\cal M}(x_{-}, x, x_{+}, J)$ of broken
connecting maps of two elements  is one. In fact, a parameterized version
of this plays an important role to prove that 
the index homology 
so defined  is independent of the choices involved and hence is an
invariant of the contact structure. In term of the parameterized
moduli space this is equivalent to show that $$ 
\dim {\tilde {\cal M}}(x_{-}, x_{+}, J)= 
\dim {\tilde {\cal M}}(x_{-}, x, x_{+}, J)-1.$$

Since the virtual dimensions of the above two  parameterized moduli
spaces can be calculated by the index formula for the corresponding
linearized ${\bar {\partial }}$-operators, the desired relation above
is equivalent to say that the index for the first moduli space
 is equal to the index  for the second moduli space minus 1.

On the other hand, in the usual Bott-type Floer homology in
symplectic case, to prove that the boundary map defined there really gives
rise to a homology theory, one needs to consider a similar relation for index
formula there. However, in that case, unlike here,  
the indices of the above two
moduli spaces are the same.

As we mentioned in [L1], the index formula above together with the
"hidden" symmetries in the target discovered in [L1], implies that
bubbling is a codimension two phenomenon and splitting of connecting
maps into broken ones of two elements is of codimension one. This will
promote various  interesting constructions in contact geometry. On the
hand, if the index formula is the same as the usual one in the Bott-type
Floer homology in symplectic case, we will have a rather different
picture. In that case, the splitting of connecting maps will become
a co-dimension two phenomenon. This will lead to a quite different
alternative for applications. For example, instead of having the index
homology outlined in [L1] and various related constructions there, one will
have a Floer homology theory and related multiplicative structure in the 
symplectization by using perturbed  connecting maps for some suitable
chosen Hamiltonian function on the symplectization.
One of the purposes of this paper is to show that the first alternative
of the consequences of [L1] can be realized.

The second motivation of this paper is to discuss a special phenomenon
concerning the dual features of the transversality of  the trivial 
connecting map. Given a closed orbit $x$, let $x_{-}=x_{+}=x$, and
consider the moduli space $ {\tilde {\cal M}}(x_{-}, x_{+}, J)$.
In this case, there is only one unparameterized connecting map, which is
the trivial map coming from $x$. It is well-known that 
$ {\tilde {\cal M}}(x_{-}, x_{+}, J)$ can be realized as  the zero set
of the 
${\bar{\partial}}$-section for some Banach bundle over a suitable
chosen Banach manifold ${\bf B}$ containing
$ {\tilde {\cal M}}(x_{-}, x_{+}, J)$. 

In order to prove the invariance of the index homology outlined above,
it is crucial to prove that the ${\bar{\partial}}$-section at the
trivial map is transversal for the  Banach manifold ${\bf B}$ naturally
appeared in the situation here. In the Floer homology in symplectic
geometry, the corresponding statement looks almost trivial. However, it
is a corner stone  of the proof of the invariance of Floer homology
with respect to the "continuation" of parameters involved in the
construction. The reason for this is that it is  the only "existence"
type statement in the proof.

We have mentioned before that in the simplest case of bubbling
off only one bubble, the domain of the limit map
will split into three components and target into two. It may happen
that the "new" principal component of the limit map is a trivial
component map. In order to prove that the virtual codimension of the
corresponding moduli space, in which  the limit map lies, is two,
it is crucial to show that the ${\bar{\partial}}$-section at the
trivial map of the principal component  is not transversal for 
a suitable choice of  Banach manifold ${\bf B}$. 

Therefore, unlike the Floer homology in symplectic geometry, the
transversality of the trivial connecting map has  double meanings here.
The second main purpose of this paper is to discuss this special
phenomenon. 

This paper is organized as follows. 

In Section. 2, we will derive the the formula for the index of the
linearized ${\bar{\partial}}$-operator and prove the transversality
and the non-transversality of the trivial connecting map.

In Section. 3, we will prove the additivity of the index of the 
${\bar{\partial}}$-operator under the deformation (gluing)
of pseudo-holomorphic maps.

\noindent {\bf Acknowledgment}: The author is very grateful
to Professor G. Tian for  valuable  and inspiring
discussions, for his help on various aspects of the project and
 for his 
encouragement.

\section{Index of linearized ${\bar{\partial}}$-operator}

As this paper is the continuation of [L1], we will assume that
the readers are  already familiar with [L1]. We now  briefly recall
some of the notations used in [L1].

Let  $(M^{2n+1}, \xi)$ be a compact contact manifold, where $\xi$ is a 
generic $2n$-dimensional subbundle of $TM$. A contact form
$\lambda=\lambda_\xi$ associated to $\xi$ is a 1-form such that
$\lambda\wedge(d\lambda)^n \not = 0$ and $\xi = ker\lambda.$
 The 2-form
$d\lambda$ is non-degenerate when restricted to $\xi$ and has a 
1-dimensional kernel at each tangent space of $M$.
  It has a canonic section
$X_\lambda$ defined by requiring that $\lambda(X_\lambda) =1$.
The vector field $X_\lambda$ is called contact field or Reeb field.
There is a decomposition $TM=\xi\oplus\bf{R}X_\lambda.$

The symplectization of $(M^{2n+1}, \xi,\lambda)$ is defined as follows. 

Let ${\widetilde M}$ be $ M\times {\bf R}$ equipped with the exact symplectic
form $\omega=d(e^r \cdot\lambda), $ where $r$ is the coordinate for
the ${\bf R}$-factor. Since $d\lambda$ is symplectic along $\xi$, there
exists a $d\lambda$-compatible almost complex structure $J$ defined on
$\xi$. In fact, the set of all such $J$'s is contractible. We extend $J$
to an $r$-invariant almost complex structure ${\tilde J}$ by
requiring:
$${\tilde J}(\frac{\partial}{\partial{r}})=X_\lambda, \, {\tilde
J}(X_\lambda)=-\frac{\partial}{\partial r}, \mbox{ and } {\tilde
J}=J$$ along $\xi.$

\vspace{3mm}

\noindent $\bullet$
{\bf Equation for ${\tilde J}$-holomorphic curves in ${\widetilde M}$}

\vspace{3mm}

Let ${\tilde u}=(u, a):\Sigma=S^1\times{\bf R}\rightarrow{\widetilde M}$
be a ${\tilde J}$-holomorphic map where $u:\Sigma\rightarrow M$ and
$a:\Sigma\rightarrow{\bf R}$. Then we have 
$$
{\tilde J}({\tilde u})\circ d{\tilde u} = d{\tilde u}\circ
i,\label{star}\eqno(\star)
$$
where $i$ is the standard complex structure on $\Sigma, $ i.e.
$i(\frac{\partial}{\partial s})=\frac{\partial}{\partial t}$ and
$i(\frac{\partial}{\partial t})=-\frac{\partial}{\partial s}.$ Here $(s,
t)$ is the cylindrical coordinate of ${\bf R}\times S^1$. 

Equation ($\star$) is equivalent to the following equations:
$$
\left\{
\begin{array}{cr}
  \pi(u) du+J(u)\pi (u) du\circ i = 0 & \quad \quad\quad(1) \\
 (u^*\lambda)\circ i=da & \quad\quad \quad(2) 
\end{array}
\right.
$$
Equation (1) is equivalent to :
$$\pi(u)(\frac{\partial u}{\partial s}) + J(u)\pi(u)(\frac{\partial
u}{\partial t})=0.\quad \quad\quad \quad(1')
$$
 
\noindent $\bullet$ {\bf Banach Space Set-up}:

\noindent
 Local coordinate near $x$:
As in [L1], we assume throughout this paper that the constant 1-form
$\lambda$ is generic in the sense that 1 is not an eigen-value of the
Pincare returning maps along the normal direction of all closed orbits
of $\lambda$. This implies that the set of unparameterized closed orbits
is discrete. 

Given a parameterized closed orbit $x:S^1={\bf R}/{\bf Z}\rightarrow M$ of
$\lambda$-period $c=\int_{S^1}(x^*\lambda)>0$, we have
$\frac{dx}{dt}=c\cdot X_\lambda(x).$ Here $X_\lambda$ is the contact field of
$\lambda$ and $t\in [0,1)={\bf R}/{\bf Z}.$ 

Let $\tau$ be the smallest positive number such that $x(\tau +t)=x(t)$,
for any $t\in [0,1).$ Let $T=\tau\cdot c.$ For any $z=x(t), t\in
[0, \tau)$, on the image of the closed orbit  $x$, we define the
$\theta$-coordinate of $z$ to be $\theta=c\cdot t, \theta\in [0, T).$
Clearly, along the image of $x$, $\frac{\partial}{\partial
\theta}=X_\lambda.$

Let

$$V_\epsilon=\{(\theta, y)\, | \, \theta\in {\bf R}/T{\bf Z}=[0, T), \, y\in
{\bf R}^{2n}, \, |y|<\epsilon\}.$$ 

Define the $(\theta, y)$-coordinate near $x$ by using the exponential
map: 
$$Exp:\, (\theta, y)\rightarrow exp_{x(\frac{\theta}{c})}(\Sigma
y_ie_i).$$
Here $y=(y_1, \cdots, y_{2n})\in {\bf R}^{2n}$ and $\{e_1, \cdots, e_{2n}\}$ is a
symplectic framing for the symplectic bundle $\xi|_x.$ 

Let $U_\epsilon(x)$ be the image of $V_\epsilon$ under the exponential
map $Exp.$ 

Now given $x_{\pm}\in {\cal P}(\lambda)$ and ${\tilde u}\in {\tilde{\cal  M}}(x_-,
x_+; J)$, we write ${\tilde u}$ as ${\tilde u}=(a, u)=(a, \theta, y)$
when $|s|$ is large enough. Here ${\cal P}(\lambda)$ is the set of
closed orbits of $\lambda$ and ${\tilde{\cal M}}(x_-, x_+; J)$ is the
moduli space of connecting ${\tilde J}$-holomorphic maps connecting
$x'_-\in \{x_-\}$ and $x'_+\in \{x_+\}$. 

Recall that in [L1], Section 4, we proved that there exist $N, \delta,
C, d_{\pm, 1} $ and $d_{\pm, 2}$ such that when $|s|>N$, we have

\begin{eqnarray*}
|y(s, t)| & < & C\cdot e^{-\delta(|s|-N)}\\
|a(s, t)-(c_\pm s+d_{\pm, 1})| & < & C\cdot e^{-\delta (|s|-N)}\\
|\theta(s,t)-(c_\pm t +d_{\pm, 2})| & < & C\cdot e^{-\delta (|s|-N)}.
\end{eqnarray*}
Here $c_\pm =\int x^*_\pm \lambda.$

Given such a ${\tilde u}\in {\tilde {\cal M}}(x_-, x_+; J)$, we want to
find a Banach manifold ${\cal B}={\cal B}^p_{1, \chi}(\tilde u)$ of
certain $L^p_1$-maps as a neighborhood of ${\tilde u}$ and a Banach
bundle ${\cal L}\rightarrow {\cal B}$ such that ${\bp}$ gives
rise to a section of the bundle and ${\tilde{\cal M}}(x_-, x_+;
J)\cap{\cal B}$ appears as a zero locus of the ${\bp}$-section.

Given $\epsilon >0$, let $$D_\epsilon(d)=\{(d_{\pm, 1}', d_{\pm, 2}')\,
|\, (d_{\pm, 1}', d_{\pm, 2}')\in {\bf R}^4\, \, , |d_{\pm, i}'-d_{\pm,
i}|<\epsilon\}.$$
For any $d'=(d_{\pm, 1}', d_{\pm, 2}')\in D_\epsilon(d), N' >N$, we
define ${\tilde u}_{d', N'}:{\bf R}^1\times S^1\rightarrow {\tilde{\cal
M}}$ connecting $x_-'=c_- t+d_{-, 2}'\in \{x_-\}$ and $x_+'=c_+t+d_{+,
2}'\in \{x_+\}$ as follows. 
Write ${\tilde u}_{d', N'}=(a, \theta, y)$, when $|s|>N.$ 

(i) When $|s|>N'+1$, define $\udn (s, t)=(c_\pm s+d_{\pm, 1}', c_\pm t
+d_{\pm, 2}', 0).$

(ii) When $|s|< N'$, define $\udn ={\tilde u}.$

(iii) When $N' < |s| < N'+1$, first write the trivial map as
$exp_{\tilde u} \xi$, and then
define $\udn(s, t)=exp_{\tilde u} \beta (s)\xi (s,t ), $ 
where $\beta (s)$ is a 
cut-off function such that $|s|<N'$, $\beta (s)\equiv 0$ and $|s|> N' +1,
\beta (s)\equiv 1.$

Let $\chi$ be the collection of the parameters, $N', \epsilon>0,
\delta>0$, we define the Banach neighborhood ${\cal B}={\cal B}^p_{1,
\chi}(\tilde u)$ of ${\tilde u}$ as follows. 
$${\cal B}=\cup_{d'\in{\cal B}_\epsilon(d)}\{{\tilde v}={\tilde
v}_{d'}=exp_{\udn}\xi | \, \xi\in L^p_{1, loc}(\udn ^*T{\tilde M}), \, 
\|\xi\|_{1, p, \delta}<\infty\}.$$ Here the weight Sobolev norm is
defined to be $\|\xi\|_{1, p,\delta}=\|e^{\delta|s|}\cdot \xi\|_{1, p}.$
Hence ${\cal B}$ is a Banach bundle over $D_\epsilon(d)$. Note that here
we assume that the exponential map $exp_{\udn}\xi$ above is the normal
exponential when $|s|$ is large enough. Hence, if we write
$$\xi=\xi_a\frac{\partial}{\partial
a}+\xi_\theta\frac{\partial}{\partial\theta}+\Sigma\xi_i\frac{\partial}{\partial
y_i}$$ in terms of the coordinate vector fields, then the exponential
map here is just the identity map along $\frac{\partial}{\partial y}$-
direction and is the identity map shifted by a constant along
${\dth}$-direction under the obvious
interpretation. This assumption will simplify some of our calculations.
When $\epsilon$ is small enough, let $\pi_{d'}$ be the parallel
transport from $\udn ^*T{\tilde M}$ to ${\tilde u}^*T{\tilde M}$ along the
shortest geodesic connecting any two points $\udn (s, t)$ to ${\tilde
u}(s, t)$. We can identify ${\cal B}$ with $L^p_{1, \delta}({\tilde
u}^*_{d, N'}T{\tilde M})\times D_\delta(d)$ by sending $(\xi , d')\in
L^p_{1,\delta}({\tilde u}^*_{d, N'}TM)\times D_\delta(d)$ to ${\tilde
V}=exp_{\udn}(\pi_{d'}(\xi))\in {\cal B}.$ 

Note that we may assume that the coordinate vector field
$(\frac{\partial}{\partial t}, \frac{\partial}{\partial y})$ is already
parallel along the two trivial maps at two ends $x_\pm.$ To simplify our
presentation, we will assume further that the almost complex structure
$J$ on $\xi$ is given by the standard complex structure $J_0$ in terms
of the framing $\{e_1, \cdots, e_{2n}\}$ on $U_\epsilon(x).$ Here
$e_i\in \xi$ such that $\pi_y(e_i)={\dyi}$, $i=1,
\cdots, 2n, $ where $$\pi_y: TM|_{U_\epsilon(x)}={\bf
R}\{{\dth}\}\oplus{\bf
R}\{\frac{\partial}{\partial y_1}, \cdots, \frac{\partial}{\partial
y_{2n}}\}\rightarrow {\bf R}\{\frac{\partial}{\partial y_1}, \cdots,
\frac{\partial}{\partial y_n}\}$$ be the bundle projection. Therefore, we
have $J(e_i)=e_{i+n}$ and $J(e_{i+n})=-e_i, 1\le i \le n.$

We define the Banach bundle ${\cal L}\rightarrow {\cal B}$ as follows.
For any ${\tilde v}={\tilde v}_{d'}=exp_{\udn} \xi\in {\cal B}$, we
define ${\cal L}_{\tilde v}=L^p_{0, \delta}({\tilde v}^*T{\tilde M})$.
Then it is clear that ${\bp}:{\cal B}\rightarrow{\cal L}$ is a
smooth section of the bundle. 

We now define the local trivialization of the bundle ${\cal
L}\rightarrow{\cal B}.$ Let ${\it p}:\udn^*T{\tilde M}\rightarrow{\tilde
v}_{d', N'}^*T{\tilde M}$ be the parallel transport of the bundle
induced from the parallel transport along geodesics from $\udn(s, t)$ to
${\tilde v}_{d', N'}(s, t).$ Here ${\tilde v}_{d', N'}^*=exp_{\udn}\xi$
as above. We define the trivialization of ${\cal L}$ by sending $(\xi,
d', \eta)\in L^p_{1, \delta}({\tilde u}^*_{d, N'}T{\tilde M})\times
D_\epsilon(d')\times L^p_{0,\delta}({\tilde u}^*_{d, N'}T{\tilde M})$ to 
${\it
p}(\pi_{d'}(\eta))\in {\cal L}_{{\tilde v}_{d', N'}}.$ Note that 
$\pi_{d'}:{\tilde u}^*_{d, N'}T{\tilde M}\rightarrow\udn T{\tilde M}$ and 
${\tilde v}_{d', N'}=exp_{\udn}\xi.$

$\bullet$ {\bf The ${\bp}$-section in local coordinate and
local trivialization:}

In terms of the local coordinate and local trivialization above,
${\bp}$-section becomes a function:
$$F:L^p_{1, \delta}({\tilde u}^*_{d, N'}T{\tilde M})\times D_\epsilon
(d)\rightarrow L^p_{0, \delta}({\tilde u}^*_{d, N'}T{\tilde M}).$$
More precisely, for any 
$$(\xi, d')\in L^p_{0, \delta}({\tilde
u}^*TM)\times D_\epsilon(d),\, \, \quad F(\xi)=({\it p}\circ\pi_{d'})^{-
1}({\bp}(exp_{\udn}\pi_{d'}(\xi))).$$
For the purpose of this section, it is important to know the formula of
$F$ in terms of the local framing defined on $U_\epsilon(x_\pm)$. We
will assume that the local framing $\{e_1, \cdots, e_{en}\}$ for
$\xi|_{U_\epsilon(\pm x)}$ is parallel and use it to define the local
trivialization of ${\cal L}$ above. Note that $e_i|_{\pm
x}={\dyi}|_{\pm x}. $ 
Then $\{{\da}, X_\lambda, e_1, 
\cdots, e_{2n}\}$ forms a local  parallel framing on $U_\epsilon(\pm
x)$.  In $U_{\epsilon}(\pm x)$,  for any given $\eta\in L^p_{o, \delta}
{\tilde u}^*_{d, N'}T{\tilde M}, $ we write
$\eta =\eta |_{U_{\epsilon}(\pm x)}=\eta_a{\da}+\eta_0
X_\lambda +\Sigma^{2n}_{i=1}\eta_i e_i.)$ 

Given $(\xi, d')\in L^p_{1, \delta}({\tilde u}^*_{d, N'}T{\tilde M})\times 
D_\epsilon(d)$, write $\xi=\xi_a\frac{\partial}{\partial
a}+\xi_0\frac{\partial}{\partial\theta}+\Sigma^{2n}_{i=1}\xi_i\frac{\partial}{\partial
y_i}$, where ${\dth}$ and
${\dyi}$ are the coordinate vector field
restricted to $x_\pm.$

Write $exp_{\udn}(\xi)=(a, \theta, y)$, then by our assumption above,
$exp_{\udn}(\xi)=(\xi_a+(c_\pm s+d_{\pm, 1}'), \xi_0+(c_{\pm}
t+d_{\pm, 2}'), (\xi_1, \cdots, \xi_{2n})).$ 

To find $F(\xi)$, we only need to write ${\bp}exp_{\udn}(\xi)$
in terms of the parallel framing $\{{\da},
X_\lambda, e_i\}$. Now 
$$\frac{\partial}{\partial s}(exp_{\udn}(\xi))=(c_\pm
+(\xi_a)_s)\frac{\partial}{\partial
a}+(\xi_0)_s\frac{\partial}{\partial
\theta}+\Sigma^{2n}_{i=1}(\xi_i)_s{\dyi}.$$
Denote $\xi_0\frac{\partial}{\partial
\theta}+\Sigma^{2n}_{i=1}\xi_i{\dyi}$ by $\eta$,
then $\eta_s=\lambda(\eta_s)X_\lambda+\pi(\eta_s), $ where $\pi$ is the
projection of $TM={\bf R}X_\lambda\oplus\gamma\rightarrow\gamma$, given by
$\pi(\eta)=\eta-\lambda(\eta)X_{\lambda_0},$
where $\gamma$ is the contact structure. It is proved in [L1] that
if we write 
\begin{eqnarray*}
\pi(\eta_s) & = & \eta_s-
\lambda(\eta_s)X_\lambda \\
 & = & \Sigma^{2n}_{i=1}((\eta_s)_i-
\lambda(\eta_s)X_i){\dyi}+((\eta_s)_0-
\lambda(\eta_s)X_0){\dth},
\end{eqnarray*}
 then 
$$\pi(\eta_s)=\Sigma^{2n}_{i=1}((\eta_s)_i-
\lambda(\eta_s)X_i)e_i=\Sigma_i((\xi_i)_s-\lambda(\eta_s)X_i)e_i.$$
Hence 
\begin{eqnarray*}
\frac{\partial}{\partial s}(exp_{\udn}(\xi)) & = & (c_\pm
+(\xi_a)_s){\da}
 + (
(\xi_0)_s\lambda(\frac{\partial}{\partial\theta})+
\Sigma^{2n}_{i=1}(\xi_i)_s\lambda(\frac{\partial}{\partial
y_i}))X_\lambda\\
&  + & \Sigma^{2n}_{i=1}((\xi_i)_s-\lambda(\eta_s)X_i) e_i.
\end{eqnarray*}
Similarly, 
\begin{eqnarray*}
\frac{\partial}{\partial t}(exp_{\udn}(\xi)) & =& (\xi_a)_t\frac
{\partial}{\partial a}+(c_\pm
+(\xi_0)_t)\frac{\partial}{\partial\theta}\\
& + & \Sigma_{i=1}^{2n}(\xi_i)_t{\dyi}
\end{eqnarray*}
and 
\begin{eqnarray*}
\frac{\partial}{\partial t}(exp_{\udn}(\xi)) & =&
(\xi_a)_t{\da}
 + ((c_\pm +(\xi_0)_t)\lambda({\dth})
  + \Sigma^{2n}_{i=1}((\xi_i)_t\lambda(\frac{\partial}{\partial
y_i}))X_\lambda\\
& + & \Sigma^{2n}_{i=1}((\xi_i)_t-(c_\pm
\lambda({\dth})+ \lambda(\eta_t))X_i)e_i.
\end{eqnarray*}
Recall that ${\tilde J}(e_i)=e_{i+n}, {\tilde J}(e_{i+n})=-e_i,$ ${\tilde
J}({\da})=X_\lambda$ and ${\tilde J}(X_\lambda)=-
{\da}.$ 

Hence, 
\begin{eqnarray*}
{\bar \partial}(exp_{\udn}(\xi)) & = & (\frac{\partial}{\partial s}+{\tilde
J}\frac{\partial}{\partial t})(exp_{\udn}((\xi))\\
& = & ((\xi_a)_s-(\xi_0)_t\lambda(\frac{\partial}{\partial
\theta}+((\xi_0)_s\lambda(\frac{\partial}{\partial\theta})+(\xi_a)_t)X_\lambda\\
& + & \Sigma^{2n}_{i=1}((\xi_i)_s-\lambda(\eta_s)X_i)e_i\\
& + & \Sigma^{2n}_{i=1}J_0\{(\xi_i)_t-(c_\pm \lambda(\frac{\partial}{\partial
\theta})+\lambda(\eta_t))X_i\} e_i\\
& - & (\Sigma^{2n}_{i=1}(\xi_i)_t\lambda(\frac{\partial}{\partial
y_i}))\frac{\partial}{\partial
a}+(\Sigma^{2n}_{i=1}(\xi_i)_s\lambda(\frac{\partial}{\partial
y_0}))X_\lambda.
\end{eqnarray*}
Therefore, in terms of $({\da},
{\dth}, {\dyi})$ along
${\tilde u}$, 
\begin{eqnarray*}
F(\xi, d') & = & ((\xi_a)_s-(\xi_0)_t\lambda(\frac{\partial}{\partial
\theta})){\da} + 
((\xi_0)_s\lambda(\frac{\partial}{\partial\theta})+(\xi_a)_t)\frac{\partial}{\partial
\theta}\\
& -& (\Sigma^{2n}_{i=1}(\xi_i)_t\lambda(\frac{\partial}{\partial
y_i}))\frac{\partial }{\partial
a}+(\Sigma^{2n}_{i=1}(\xi_i)_s\lambda(\frac{\partial}{\partial
y_i})){\dth}\\
& + & (\Sigma^{2n}_{i=1}((\xi_i)_s-\lambda(\eta_s)X_i)\dyi\\
& + & \Sigma^{2n}_{i=1}J_0\{(\xi_i)_t-
(\lambda(\eta_t)+c_\pm\lambda(\dth))X_i\}\dyi.
\end{eqnarray*}

Note that here $\lambda(\dth)$, $\lambda(\dyi)$, $\lambda(\eta_s)$ and
$\lambda(\eta_t)$ are evaluated along $exp_{\udn}\xi.$ 

\vspace{3mm}

\noindent $\bullet$ {\bf  Linearization of the ${\bar {\partial }}$-operator}:

\vspace{3mm}

Now 
\begin{eqnarray*}
DF(0)(\xi, d') & = & \lim_{r\rightarrow 0}\frac{F(r\xi, d+rd')-F(0,
d)}{r}\\
& =& \lim_{r\rightarrow 0}\frac{F(r\xi, d+rd')}{r}
\end{eqnarray*}
when $|s|> N.$ 

Since for $|s|>N$, when $r\rightarrow 0$,
$exp_{{\tilde u}_{d+rd', N'}}r\xi\rightarrow{\tilde u}_{d, N'}$ of
 the trivial map near $x_\pm$, it is easy to
see that the $(\da, \dth)$-component of 
$$DF(0)(\xi, d')=((\xi_a)_s-(\xi_0)_t)\da+((\xi_0)_s+(\xi_a)_t)\dth.$$ 
To find the $y$-component of $DF(0)(\xi, d')$, we use the matrix
notation. Let 
$${\underline \xi}=\left (\begin{array}{l}
	\xi_1\\
	\vdots\\
	\xi_{2n}
	\end{array}
	\right ) \quad \quad \mbox{and}\quad\quad Y=\left (\begin{array}{l}
	X_1\\
	\vdots\\
	X_{2n}
	\end{array}
	\right ).$$
We proved in [L1] that in local $(\theta, y)$-coordinate, 
$$Y(\theta,
y)=\int^1_0 dY(\theta, \tau y)d\tau\left (\begin{array}{l}
	y_1\\
	\vdots\\
	y_{2n}
	\end{array}
	\right ),$$
where $dY(\theta, y)$ is a $2n\times 2n$ matrix whose $(i, j)$-entry is
$\frac{\partial X_i}{\partial y_j}.$ 

Denote $\int^1_0 dY(\theta, \tau y)d\tau$ by $DY(\theta, y)$. 

In our case, for $exp_{{\tilde u}_{d+rd', N'}}(r\xi)$, 
$\theta=(c_\pm t + d+rd')+\xi_0 , $ and $y_i=r\cdot\xi_i, i=1, \cdots, 2n,$ 
the $y$-component of 
$$DF(0)(\xi, d')={\underline \xi}_s+J_0{\underline
\xi}_t-c_\pm J_0 dY(c_\pm t +d_{\pm, 2}, 0){\underline \xi}.$$
In [L1], we have assumed that $c_\pm=1$ by rescaling. It was proved
there that $S=-c_\pm J_0 dY(t, 0)$ is a symmetric matrix and all the eigen
values of the self-adjoint elliptic operators $A: L^2_1(S^1, {\bf
R}^{2n})\rightarrow L_0^2(S^1; {\bf R}^{2n})$ defined by
$A(z)=J_0\frac{dz}{dt}+S\cdot z$ are non-zero when $\lambda$ is generic.

Hence the $y$-component of $DF(0)(\xi, d')={\underline \xi}_s+A\cdot
{\underline \xi}.$

Write ${\tilde \xi}=(\xi_a, \xi_0)$ and identify it with $\xi_a+\xi_0
i.$ Then the $(\da, \dth)$-component of $DF(0)(\xi,
d')$ is given by ${\bp}\cdot{\underline \xi}=(\ds+i\dt)(\xi_a+\xi_0 i)$
if we identify $\{\da, \dth\}$ with $\{1, i\}.$ 

Hence when $|s|>N$, 
$$DF(0)(\xi, d')={\tilde L}_1({\tilde \xi})+{\tilde L}_2({\underline
\xi}),$$ when $\xi=({\tilde \xi}, {\underline \xi})$, 
${\tilde L}_1={\bp}=(\ds+i\dt)$ and ${\tilde L}_2=\ds+A.$ 

For an arbitrary $s$, the formula for $DF(0)(\xi, 0)$ is well-known. We
have $DF(0)(\xi, 0)=\xi_s+{\tilde J}({\tilde u})\xi_t+B({\tilde
u})\cdot\xi,$ where ${\tilde J}=\left [\begin{array}{lr}
					i & 0\\
					0 & J	
					\end{array} \right ]$ with 
respect to the decomposition $T{\tilde M}={\bf R}\{\da\}\oplus{\bf
R}X_\lambda\oplus\gamma$, and $B(\tilde u)$ is some matrix operator. 
$$DF(0)(0, d')=\lim_{r\rightarrow 0}\frac{1}{r}\{{\it p}^{-
1}({\bp}{\tilde u}_{d+rd', N'})-{\it p}^{-
1}({\bp}{\tilde u}_{d, N'})\}.$$
Note that when $|s|>N$, $\bp({\tilde u}_{d+rd', N'})=0=\bp {\tilde u}_{d, N'}$,
$DF(0)(0, d')\equiv 0.$ In general, $DF(0)(0, d')=C({\tilde u}_{d, N'})\cdot d'$ for
some matrix operator $C.$ Therefore, $DF(0)$ is an elliptic operator
over ${\bf
R}^1\times S^1$, a particular non-compact manifold with cylindrical
ends.

In this situation, it is proved in [LM] that the elliptic operator
$DF(0):L^p_{1, \delta}({\tilde u}_{d, N'}^* T{\tilde M})\times {\bf R}^4\rightarrow
L^p_{0, \delta}({\tilde u}_{d, N'}^* T{\tilde M})$ is Fredholm if and only if
$\delta$ is not an eigenvalue of the operator $i\dt\oplus A:L^2_1
(S^1,{\bf R}^{2n+2})\rightarrow L^2_0(S^1, {\bf R}^{2n+2}).$

Write $B=\left (\begin{array}{lr}
		B_{a, X_\lambda} & B_1\\
		B_2 & B_\xi
		\end{array} \right ).$
It is easy to see that we can deform $DF(0)$ to get rig of the $B_1, B_2
$ and $B_{a, X_\lambda}$ terms and $Cd'$ in the middle part but maintain
the same asymptotic behavior of $DF(0)$. Then each operator in the
deformation is still Fredholm and therefore has the same index. 

Let $L$ be the resulting operator. Then, $L=L_1\oplus L_2.$  Hence,
$$L_1: L^p_{1, \delta}({\tilde u}_{d, N'}^*({\bf R}\{\da\}\oplus{\bf
R}X_\lambda))\oplus{\bf R}^4=L^p_{1, \delta}({\bf R}^1\times S^1, {\bf
R}^2)\oplus {\bf R}^4\rightarrow L^p_{0, \delta}({\bf R}^1\times S^1,
{\bf R}^2)$$ given by $L_1({\tilde \xi}, d)=\bp {\tilde \xi}.$ and 
$$L_2:L^p_{1,\delta}({\tilde u}_{d, N'}^*(\xi))\rightarrow L^p_{0,
\delta}({\tilde u}^*_{d, N'}(\xi))$$ given by $L_2({\underline
\xi})=(\underline \xi)_s+J({\tilde u}_{d, N'})(\underline
\xi)_t+B_\xi{\underline \xi}.$

Since $0$ is not an eigenvalue of $A$, ${\underline L}_2: L^p_1({\tilde
u}_{d, N'}^*(\xi))\rightarrow L^p_0({\tilde u}_{d, N'}^*(\xi))$ is also
Fredholm with the same index as $L_2$. We will still use $L_2$ to denote
${\underline L}_2$. Note that as far as the computation of index is
concerned, we can replace $p$-norm, by $2$-norm.

We will denote the restriction of $L_1$ to $L^p_{1, \delta}({\bf
R}^1\times S^1; {\bf R}^2)$ by $L_{1, \delta}.$ For the application of
[L2], we need to consider the multi-ends case also.

Given $l+m$ closed orbits $x_{-, 1}, \cdots, x_{-, l}$, $x_{+, 1},
\cdots, x_{+, m}$, with each $x_{\pm , j}:S^1\rightarrow M$ satisfying
the equation $\frac{d x_{\pm, j}}{dt}=c_{\pm, j}X_\lambda(x_{\pm, j})$,
we denote the collection $\{x_{-, 1}, \dots, x_{-, l})\}$ by $x_-$ and
$\{x_{+, 1}, \cdots, x_{+,m}\}$ by $x_+.$ 

Let $(S^2, z_{-1, 1}, \cdots, z_{-, l}, z_{+, 1}, \cdots, z_{+, m})\in
{\cal M}_{0, l+m}.$ We identify a small punctured disc $D_\epsilon(z_{-, i})
\setminus\{z_{-, i}\}$ of $S^2$ with ${\bf R}^{-}\times S^1$, $i=1,
\cdots, l$ and $D_\epsilon(z_{+, j})\setminus\{z_{+, j}\}$ with ${\bf
R}^+\times S^1$, $j=1, \cdots, m$, so that $(S^2; z_-, z_+)$ can be
thought as a punctured $S^2$ with $l$-negative cylindrical ends and $m$
positive ends. Here $z_-$ stands for $\{z_{-, 1}, \cdots, z_{-, j}\}$
and $z_+$ for $\{z_{+, 1}, \cdots, z_{+, m}\}$. 
As before, we define 
\begin{eqnarray*}
& & {\tilde{\cal M}}(x_-, x_+; J) = \\
& &\{{\tilde u}:S^2\setminus \{z_-,
z_+\}\rightarrow{\tilde M}, \bp_{\tilde J}{\tilde u}=0, \mbox{ and }
\lim_{y\rightarrow{z_\pm}}{\tilde u}(y)=x'_\pm\in\{x_\pm\}\}, 
\end{eqnarray*}
where $\{x_\pm\}$ is the equivalent class of $x_\pm.$ Note that here we
allow $(S^2; z_-, z_+)$ to vary in ${\cal M}_{0, l+m}.$ As before, we
have the exponential decay estimate for ${\tilde u}$ along each of its
cylindrical ends. Given $N'>N$, $\epsilon >0$, and $D_\epsilon(d)$, we
can define $\udn$, and the Banach neighborhood ${\cal B}={\cal B}^p_{1,
\chi}(\tilde u)$ of $\tilde u$ in exactly the same way as before, where
$\chi=(N', \epsilon, \delta).$ Note that $d' $ and $d\in {\bf
R}^{2(l+n)}.$ We define the Banach bundle ${\cal L}\rightarrow{\cal B}$
slightly different from the definition before. For any ${\tilde v}\in
{\cal B}$, ${\cal L}_{\tilde v}=L^p_{0, \delta}({\tilde v}^*(\wedge^{0,
1}_{(i, \tilde J)}T{\tilde M}))$. Then $\bp$-operator can be thought as
a global section of the bundle ${\cal L}\rightarrow{\cal B}$. Here
$\bp$-section is defined by 
$$\bp ({\tilde v})=d{\tilde v}+{\tilde J}(\tilde v)\circ d{\tilde
v}\circ i.$$
Arguing the same way as before, we can find a coordinate chart for
${\cal B}$ by identifying it with an open set of $L^p_{1,\delta}(u^*_{d,
N'}TM)\times D_\epsilon(d)$ and a trivialization of ${\cal L}$ by
identifying it with
$$L^p_{1,\delta}(u^*_{d, N"}T{\tilde M})\times D_\epsilon(d)\times
L^p_{0, \delta}(u^*_{d, N'}(\wedge^{0,1}_{(i, {\tilde J})}T(\tilde M))$$
and 
$$DF(0):L^p_{1,\delta}(u^*_{d, N'}T{\tilde M})\times
D_\epsilon(d)\rightarrow L^p_{0,\delta}(u^*_{d, N'}(\wedge^{0,1}_{(i,
{\tilde J})}T(\tilde M)).$$ 
Finally, we can deform the Fredholm operator $DF(0)$ into $L=L_1\oplus
L_2$ of same index where
$$L_1:L_{1,\delta}^p(S^2\setminus\{z_-, z_+\}, {\bf R}^2)\oplus {\bf
R}^{2(l+m)}\rightarrow L^p_{1,\delta}(\wedge^{0,1}(S^2\setminus \{z_-,
z_+\})$$
given by $L_1({\tilde \xi}, d')=\bp {\tilde \xi}$, and 
$$L_2: L^p_{1,\delta}(u^*_{d, N'}(\xi))\rightarrow
L^p_{0,\delta}(u^*_{d, N'}(\wedge^{0,1}_{(i, {\tilde J})}\xi))$$ 
given by $L_2(\underline \xi)=\bp_J{\underline \xi}+{\tilde
B}\cdot{\underline \xi}. $ Here $\bp{\tilde
\xi}=(\frac{\partial {\tilde \xi}}{\partial {\bar z}})d{\bar z}$ in terms of the
local complex coordinate $z$ of $S^2\setminus \{z_-, z_+\}$,  
${\bp}_J{\underline \xi}=\nabla\xi+J(u_{d, N'})\circ\nabla\xi\circ i$
 and ${\tilde B}$ is a bundle map ${\tilde u}^*_{d,, N'}(\xi)\rightarrow 
{\tilde u}_{d,
N'}^*(\wedge^{0,1}_{(i, {\tilde J})}\xi)$. Now near $x_\pm$, by using
the local coordinate $z$ of the domain near the ends, we can identify
${\tilde u}_{d,
N'}^*(\wedge^{0,1}_{(i, {\tilde J})}\xi)
=u^*_{d, N'}(\xi)\oplus{\bf C}
\{d{\bar z}\}$ with $u^*_{d, N'}(\xi)$. Under this identification, in
the local
coordinate of $U_\epsilon(x_\pm)$, 
$$L_2(\underline \xi)=(\underline
\xi)_s+J({\tilde u}_{d, N'})(\underline \xi)_t+B_{\underline
\xi}{\underline \xi}$$ as before.

Now the virtual dimension of ${\tilde{\cal M}}(x_-, x_+; J, z_-, z_+)$ 
is equal to the index of  $L$, which is the same as $Ind(L_1)+Ind(L_2).$
Here ${\tilde {\cal M}}(x_-, x_+; J, z_-, z_+)$ is the subspace of
${\tilde {\cal M}}(x_-, x_+; J)$ whose elements have the fixed domain
$S^2\setminus\{z_-, z_+\}. $ Our goal is to study the behavior of
$Ind(L)$ under deformation (gluing) of the domain of $\tilde u$. To this
end, we will calculate $Ind(L_1)$ first. 

Recall $$L_1=\bp: L^p_{1, \delta}(\sz, {\bf C})\oplus{\bf
R}^{2(l+m)}\rightarrow L^p_{0,\delta}(\wedge(\sz))$$
and $L_{1,\delta}=L_1$ restricted to $L^p_{1,\delta}(\sz, {\bf C})$.
Note that $0<\delta<2\pi$ of the firsts eigenvalue of the operator
$i\dt$ along each ends. Here the weighted Sobelev space
$L^p_{1,\delta}(\sz, {\bf C})$ and $L^p_{0,\delta}(\wedge^{0,1}(\sz))$
are measured with respect to the cylindrical metric along each end of
$\sz .$

\begin{lemma}
If $\# \{z_-, z_+\}=l+m$, then $Ind(L_{1,\delta})=2-2(l+m).$ Therefore,
$Ind(L_1)=2. $
\end{lemma}
\proof
Start with the case $l=m=1.$ Then $\sz={\bf R}^1\times S^1$ with
cylindrical coordinate $z=(s, t)$. In terms of the global complex
coordinate $z$, $d{\bar z}$ is a global section of $\wedge^{0,1}({\bf
R}^1\times S^1)$. Hence $\wedge^{0,1}({\bf R}^1\times S^1)\equiv ({\bf
R}^1\times S^1)\times {\bf C}$ of the trivial bundle and under this
identification $L_{1, \delta}:L^p_{1,\delta}({\bf R}^1\times S^1, {\bf
C})\rightarrow L^p_{0,\delta}(\wedge^{0,1}({\bf R}^1\times S^1))$
becomes $\bp: L^p_{1,\delta}({\bf R}^1\times S^1; {\bf C})\rightarrow
L^p_{0,\delta}({\bf R}^1\times S^1; {\bf C})$ given by $\xi\rightarrow
\frac{\partial \xi}{\partial s}+i\frac{\partial\xi}{\partial t}.$ 

We have defined $L^p_{1,\delta}({\bf R}^1\times S^1, {\bf C})$ by
requiring that $\|\xi(s, t)\cdot e^{\delta |s|}\|_{1, p}<\infty$ along
both of its ends. 
We define $L^p_{1, (-\delta, \delta)}({\bf R}^1\times S^1, {\bf C})$ by
requiring $\|\xi(s, t)\cdot e^{\delta s}\|_{1, p}<\infty$ for any of its
 element $\xi$.
This is the
same as requiring $\|\xi(s, t)\cdot e^{\delta |s|}\|_{1, p}<\infty$
along the positive end $z_+$ and $\|\xi(s, t)\cdot e^{-\delta
|s|}\|_{1,p}<\infty$ along the negative end. In general, we can consider
multi-ends with arbitrary choices of the sign before $\delta$ and define
$L^p_{1, (\pm\delta_{z_-}, \pm\delta_{z_+})}(\sz, {\bf C})$ 
and $L^p_{0, (\pm\delta_{z_-}, \pm\delta_{z_+})}(\wedge^{0,1}(\sz)).$
Let $L_{1, (\pm\delta_{z_-}, \pm\delta_{z_+})}$ be the corresponding
operator.

\begin{sublemma}
$Ind(L_{1, (-\delta, \delta)})=0.$
\end{sublemma}

\proof
Define the isomorphism:
$$e_0:L^p_0({\bf R}^1\times S^1, {\bf C})\rightarrow L^p_{0, (-\delta,
\delta)}({\bf R}^1\times S^1, {\bf C})$$ and
$$e_1: L^p_1({\bf R}^1\times S^1, {\bf C})\rightarrow L^p_{1, (-
\delta, \delta)}({\bf R}^1\times S^1, {\bf C})$$ by sending
$\xi\rightarrow e^{-\delta s}\xi.$ 

Consider 
$$
D=e^{-1}_0\circ L_{1, (-\delta, \delta)}\circ e_1= e_0^{-1}\circ\bp\circ 
e_1:  L^p_1(({\bf R}^1\times S^1, {\bf C})\rightarrow L^p_0({\bf
R}^1\times S^1, {\bf C}).
$$
For any  $\xi\in L^p_1({\bf R}^1\times S^1, {\bf C})$, 
$$
D(\xi)=e^{\delta s}\bp (e^{-\delta s}\xi)=\bp\xi-\delta
\xi=\ds\xi+(i\dt-\delta)\xi.
$$
Clearly $Ind(L_{1, (-\delta, \delta)})=Ind (D).$ 

Now consider the self-adjoint operator
$$i\dt+\delta: L^2_1(S^1; {\bf C})\rightarrow L^2_0(S^1; {\bf C}).$$
Since $0<\delta <2\pi$, $0$ is not an eigenvalue of $i\dt-\delta.$ Write 
$D=\ds+L_t$, where $L_t=i\dt-\delta.$ Define the resolvent
$R_\lambda=(L_t-\lambda i)^{-1}:L^2_1(S^1; {\bf C})\rightarrow
L^2_1(S^1; {\bf C})$ for any $\lambda\in {\bf R}$, which is bounded
independent of $\lambda.$ It is proved in [F] and [LM] that in the
situation,  if we define ${\hat\xi}(\lambda)=\int e^{i\lambda s}\xi(s)ds$
and $\eta(s)=\frac{1}{2\pi}\int e^{-i\lambda s}R_\lambda
{\hat\xi}(\lambda)d\lambda.$
Then $D\eta=\xi.$ That is $D$ is invertible. Hence $Ind(D)=0$ with
respect to $L^2$-norm. Same is true for $L^p$-norm (see [LM]).

It is proved in [LM] that $Ind(L_{1, (-\delta, \delta)})-Ind(L_{1,
(\delta, \delta)})=i(-\delta_{z_-}, \delta_{z_-})$, 
where $i(-\delta_{z_-}, \delta_{z_-})$ is the dimension of the eigen
spaces of $i\dt: L^2_1(S^1, {\bf C})\rightarrow L^2_1(S^2, {\bf C})$
of the asymptotic operator associated to the end $z_{-}$
with eigenvalue in $(-\delta, \delta)$. It is the same as the dimension
of the space of solutions $\xi$ of $D$ of the form 
$\xi(s,t)=exp(\lambda s)\cdot p(s,t)$ with $-\delta <\lambda<\delta$ and
$p(s,t)$ is a polynomial in $s$ with coefficient in $C^\infty(S^1; {\bf
C})$ (See [LM]). This implies that $Ind(L_{1, (\delta, \delta)})=-2$ as 
$i(-\delta_{z_-}, \delta_{z_-})=2.$ 

Now we consider the general multi-ends case. Consider first the case 
$$
L_{1,
(-\delta, -\delta)}:
L^2_{1, (-\delta_{z_-},-\delta_{z_+})}(\sz, {\bf C})\rightarrow L^2_{0,
(-\delta_{z_-}, -\delta_{z_+})}(\wedge^{0,1}(\sz)).
$$
Let $(s, t)$ be the cylindrical coordinate near an end $p\in\{z_-,
z_+\}$ and $w\in {\bf C}$ the complex coordinate of $S^2$ near $p$ with
$w=0$ at $p$. Then $w=e^{-2\pi (s+it)}$. 
Given $\xi\in L^2_{1, (-\delta_{z_-}, -\delta_{z_+})}(\sz, {\bf C}),$
we have
$$
{\int\int}_{{\bf R}^+\times S^1}|e^{-\delta\cdot s}\xi|^2 dsdt<0.
$$
Write $w=x+iy$, $z=s+it$. Then 
$$2dx\wedge dy=idw\cdot d{\bar w}=4\pi^2\cdot i e^{-4\pi s}dz\wedge
d{\bar z}=4\pi^2 e^{-4\pi s}\cdot 2 ds\wedge dt.
$$
Then 
\begin{eqnarray*}
{\int\int}_{{\bf R}^+\times S^1}|e^{\delta \cdot s}\xi|^2 ds\wedge dt 
& = & \frac{1}{4\pi^2}\int_{D_1\setminus\{0\}}|e^{(2\pi-\delta)\cdot
s}\cdot \xi|^2 dx\wedge dy\\
& = &
\frac{1}{4\pi^2}\int_{D_1\setminus\{0\}}|\frac{\xi}{|w|^{\frac{2\pi-
\delta}{2\pi}}}|^2 dx\wedge dy
\end{eqnarray*}
 
This implies that $\int_{D_1\setminus\{0\}}|\xi|^2 dx\wedge dy<\infty.$

Assume that $\xi\in ker(L_{1, (-\delta, -\delta)})$, then $\xi$ is
holomorphic on $\sz$. The elliptic regularity implies that $\xi$
extends over $\{z_-, z_+\}$ and is well-defined on $S^2.$ From this, one
can easily conclude that $dim (ker (L_{1, (-\delta, -\delta)}))=2. $

Let $L^*_{1, (\delta ,\delta)}$ be the dual of $L_{1, (-\delta, -
\delta)}.$ Then 
$$
L^*_{1, (\delta, \delta)}: L^2_{1, (\delta,
\delta)}(\wedge^{0,1}(\sz))\rightarrow L^2_{0, (\delta, \delta)}(\sz,
{\bf C}).
$$
$Dim (ker ( L^*_{1, (\delta, \delta)}))=dim (coker L_{1, (-\delta, -
\delta)})$. As before, for any $p\in \{z_-, z_+\}$, let $z=s+it$ and
$w=e^{-2\pi z}$ be the two types of coordinates. Given 
$\eta\in  L^2_{1, (\delta, \delta)}(\wedge^{0,1}(\sz))$ near $p$, in
terms of the local section $d{\bar z}$,  $\eta=\phi\cdot d{\bar z}=\psi d {\bar w}$ with $\phi\in
L^2_1({\bf R}^+\times S^1; {\bf C})$ and $\psi\in
L^2_1(D_1\setminus\{0\}, {\bf C}).$ 
Then $\phi=\psi\cdot e^{-2\pi {\bar z}}.$

Hence 
\begin{eqnarray*}
\|\eta|_{{\bf R}^+\times S^1}\|^2_{2, \delta} & = & \int_{{\bf
R}^+\times S^1} |e^{\delta s}\eta|^2 ds\wedge dt\\
& = & \int_{{\bf R}^+\times S^1} |e^{\delta s}\phi|^2 ds\wedge dt\\
& = & \int_{{\bf R}^+\times S^1} |e^{-(2\pi-\delta)s}\cdot \psi|^2
ds\wedge dt\\
& = & \frac{1}{4\pi^2} {\int}_{D_1\setminus\{0\}} e^{4\pi s}\cdot e^{-
2(2\pi-\delta))s}|\psi|^2 dx\wedge dy\\
& = & \frac{1}{4\pi^2}{\int}_{D_1\setminus\{0\}}
|\frac{\psi}{|w|^{\frac{\delta}{2\pi}}}|^2 dx\wedge dy<\infty.
\end{eqnarray*}
Again this implies that $\|\eta|_{D_1\setminus\{0\}}\|_2<\infty.$

We want to show that if $L^*\eta=0$, then $\eta$ extends over to $p$ and
there is no restriction on the value of $\eta$ at $p$. 

To see the second statement, assume that $\eta=\psi d{\bar w}$, with
$\psi\in C^\infty(D_1, {\bf C}). $ Then we have shown that in terms of
cylindrical coordinate $z, \eta=\phi d{\bar z}$ with $\phi=\psi\cdot
e^{-2\pi(s-i t)}. $ Clearly, for $0<\delta <2\pi$, 

\begin{eqnarray*}
\|\eta|_{{\bf R}^1\times S^1}\|^2_{2,{\delta}} 
& = & \int_{{\bf R}^1\times S^1} |
e^{(\delta-2\pi)\cdot s}\psi|^2 ds\wedge dt \\
& \leq  & C \int_{{\bf R}^1\times S^1} e^{2(\delta-2\pi)\delta} ds\wedge
dt <\infty.\\  
\end{eqnarray*}
To see the first statement, we write down the explicit expression for
$L^*.$ As in [LM], we use the cylindrical Hermitian metric on $(\sz, i)$
to introduce a Hermitian inner product on  the bundle $\bigwedge^{p,
q}(\sz)$, which gives a hermitian product $<, >$ on 
$\Omega^{p,q}_0(\sz)$ of smooth sections of $\bigwedge^{p,q}(\sz)$ 
with compact support.
Now $L=L_{1, (-\delta, -\delta)}$ is just the $\bp$-operator:
$$
C^\infty_0(\sz, {\bf
C})=\Omega^{0,0}_0(\sz)\rightarrow\Omega^{0,1}_0(\sz).
$$
Similarly to the compact case (see [GH] for example), with respect to
pairing $<>$, 
$$
L^*=\bp^*=*_z\cdot \bp \cdot *_z:\Omega^{0,1}_0(\sz)\rightarrow
C^\infty_0(\sz, {\bf C}),$$
 where $*_z$ is the Hodge-star operator. Here we use
the subscript $z$ to indicate that the $*$-operator is taken with
respect to the cylindrical metric. Since the cylindrical metric is
conformal to the standard metric on $S^2$ and $*$ is conformal invariant
on $\bigwedge^{0,1}$, we have $*_z=*_w$ of the $*$-operator with respect
to the standard metric of $S^2$ on $\bigwedge^{0,1}$. On the other hand,
on the bundle $\bigwedge^{1,1}$ of top forms, we have 
$*_z=e^{-4\pi s}\cdot * .$ Hence $L^*=\bp^*_z=e^{-4\pi s}\cdot
{\bp}^*_w, $ where $\bp^*_w$ is the dual of $\bp$ with respect to the
standard metric of $S^2$, i.e. the usual dual of $\bp$ in Hodge theory.
Therefore, if $\eta\in ker (L^*_{1, (\delta, \delta)})=ker (\bp^*_z)$,
then $\eta\in L^2_{1, (\delta ,\delta)}(\wedge^{0,1}(\sz))$ and
${\bp}^*_w(\eta)=0$. Since $\bp^*_w$ is elliptic and
$\|\eta|_{D_1\setminus\{0\}}\|_2<\infty$,  $\eta$ extends over
all singular points $z_-, z_+$. As we mentioned before that there is no
restriction on the values of $\eta$ at $z_-$ and $z_+$, we conclude that 
$ker ( L^*_{1, (\delta,\delta)})$ can be identified with $ker ( \bp^*_w)$, where
$\bp^*_w: L^2_1(\bigwedge^{0,1}(S^2))\rightarrow
L^2_0(\bigwedge^{0,0}(S^2)).$ Here the $L^2$-norm is measured with
respect to the standard metric of $S^2.$ By Hodge theory,
$ker(\bp^*_w)=H^{0,1}_{\bp}(S^2)\equiv H^0({\bf C}P^1, {\cal O} (-2))=0.$
 This proves that $Ind(L_{1, (-\delta, -\delta)})=2.$ Hence 
$Ind(L_{1,\delta})=2-2(l+m)$, and $Ind(L_1)=2.$ 

Finally, we deal with
the special case where $l=0, m=1.$
Again, consider $L_{1,-\delta}$ first. Clearly, if $\eta\in ker(L_{1,-
\delta})$, then $\eta$ is holomorphic defined over ${\bf C}$ and the
ratio of growth of $\eta$ in the cylindrical coordinate $z=s+it$ as 
$s\rightarrow \infty,$
 is
less than $e^{+\delta s}.$ This implies that $\eta\equiv \mbox{ constant
}.$ Hence $dim (ker(L_{1, -\delta}))=2.$ The above argument for general
case on the kernel of the dual operator is also applicable here. We have
\begin{eqnarray*}
coker(L_{1, -\delta}) & \equiv & ker(L_{1,\delta}^*)\\
& \equiv & ker(\bp^*_w)\equiv H^0({\bf C}P^1, {\cal O}(-2))\equiv 0.
\end{eqnarray*}
Hence $Ind(L_{1,\delta})=Ind(L_{1,-\delta})-2=0$ and $Ind(L_1)=2.$
\QED

\vspace{3mm}

\noindent $\bullet$ {\bf  Transversality and non-transversality of
the trivial connecting map}:

\vspace{3mm}

Now we consider the special case that  $u\in {\tilde {\cal M}}(x_-, x_+;
J)$ with $\{x_-\}=\{x_+\}.$ In this case, since $\int_{S^1}(x_-
')^*\lambda=\int_{S^1}(x_+')^*\lambda =c$ for any $x_\pm'\in \{x_\pm\}$,
the $E_\lambda$-energy 
$$E_\lambda(\tilde u)=\int_{{\bf R}^1\times S^1}{\tilde u}^*d\lambda=0,
$$
and we have ${\tilde u}(s, t)=(cs+d_1, ct+d_2, 0)$ in the $(a, \theta,
y)$ coordinate (for the proof, see,  for example, [L1]).  Therefore, the
two asymptotic limits $\lim_{s\rightarrow\pm\infty}
u(s,t)=ct+d_2=x_\pm'\in\{x_\pm\}, $ which are the same parameterized
closed orbits, and $\tilde u$ itself is the trivial map. 

We want to calculate the index of $L$, which is the linearization $DF$ 
at $\tilde u$. To this end, we only  need to know $Ind(L_2)$ in this
case. 

Recall that $L_2:L^p_{1,\delta}({\tilde u}^*\xi)\rightarrow
L^p_{0,\delta}({\tilde u}^*\xi)$ given by sending $\eta$ to
$\ds\eta+J_0\dt\eta+B\eta.$ Note that the local coordinate
$(a,\theta,y)$ is defined over the whole image of ${\tilde u}$. In terms of
the coordinate framing, ${\tilde u}^*\xi=<\frac{\partial}{\partial y_1},
\cdots, \frac{\partial}{\partial y_{2n}}>|_{\tilde u}.$

Since $\frac{\partial}{\partial y_i}|_{\tilde u}=e_i|_{\tilde u}$, in
terms of $\{\frac{\partial}{\partial y_i}\}, $
$\eta=\Sigma^{2n}_{i=1}\eta_i \frac{\partial}{\partial y_i}$, and the
$d\lambda$-compatible almost complex structure is given by the standard
complex structure on ${\bf R}^{2n}.$ $B=B(s,t)$ is a self-adjoint matrix
operator. Hence in terms of the framing $\{\frac{\partial}{\partial
y_1}, \cdots, \frac{\partial}{\partial y_{2n}}\}$, 

$$L_2:L^p_{1,\delta}({\bf R}^1\times S^1, {\bf R}^{2n})\rightarrow
L^p_{0,\delta}({\bf R}^1\times S^1, {\bf R}^{2n}).$$

Recall that $0$ is not an eigenvalue of $J_0\dt+B:L^2_1(S^1; {\bf
R}^{2n})\rightarrow L^2_0(S^1; {\bf R}^{2n}).$ It follows from [LM] that
here 
we can use the usual Sobelev norm rather than the $\delta$-weighted one.
We can also use $L^2$-norm rather than $L^p$-norm for the purpose of
calculating the index. 

\begin{lemma}

$Ind(L_2)=0$. In fact, 
$$
L_2: L^2_1({\bf R}^1\times S^1, {\bf R}^{2n})\rightarrow L^2_0({\bf
R}^1\times S^1, {\bf R}^{2n})$$
 given by $\eta\rightarrow
\ds\eta+J_0\dt\eta+B\eta$ is an isomorphism. 

\end{lemma}
\proof

Write $L_2=\ds+L_t$, where $L_t=J_0\dt+B$ is self-adjoint. As before, we
define the resolvent 
$$
R_\lambda=(L_t-i\lambda)^{-1}:L^2(S^1; {\bf R}^{2n})\bigotimes{\bf
C}\rightarrow L^2_1(S^1; {\bf R}^{2n})\bigotimes {\bf C},$$
which is bounded independent of $\lambda\in{\bf R}.$ Define
${\hat\eta}(\lambda)=\int e^{i\lambda s}\eta(s) ds.$ Then
$\gamma(s)=\frac{1}{2\pi}\int e^{-i\lambda s}
R_\lambda{\hat\eta}(\lambda)d\lambda$ has the property that
$L_2\gamma=\eta.$ It follows that $L_2$ is invertible (See [LM] and [F]
for more details).

\QED

Therefore, $Ind(L)=Ind(L_1)+Ind(L-2)=2.$ 

\medskip

Recall that we embedded a neighborhood of $\tilde u$ in ${\tilde{\cal
M}}(x_-, x_+; J)$ into ${\cal B}={\cal B}^p_{1, \chi}(\tilde u)$ and
defined ${\cal L}\rightarrow{\cal B}$ with the $\bp$-section. 
To define
$\cal B$, we introduced $\udn$ before. 
It is easy to see that in the case $\{x_-\}=\{x_+\}, $ those $\udn$'s
are very close to the trivial connecting map ${\tilde u}_d$ if $|d-d'|$
is very small. Note also that the linearized $\bp$-section at ${\tilde
u}_d$, $DF(0)$ is just $L=L_1+L_2$ in this case. We have shown that
$ker(L_2)=0$, $coker(L_2)=2$ and $Ind(L)=2.$ Intuitively it is more or less
clear that $ker
(L)=2.$  To get a precise prove, we directly calculate the $coker (L_{1})$.
By the definition of $DF(0)$, a short computation shows that
$DF(0) (0, d')=\beta_{+}' d'_{+}+\beta_{-}' d'_{-}, $ for any $d'=
(d_{+}, d_{-})\in R^4$ and $d_{\pm} \in {\bf R}^2={\bf C}$, 
where $\beta_{\pm}$ is the cut-off function we introduced before.
This proves that $coker (L_{1}) =0$, and we have the following propose
tion on the transversality of the trivial connecting maps.

\begin{lemma}

When $\{x_-\}=\{x_+\}$, $\bp:{\cal B}\rightarrow{\cal L}$ is a
transversal section and its zero locus ${\tilde{\cal M}}(x_-,x_+; J)\cap
{\cal B}$ is a smooth manifold of $2$-dimensional. Moreover
$\bp:{\underline {\cal B}}\rightarrow {\cal L}$ is also a transversal
section over the unparameterized moduli space ${\cal M}(x_-, x_+; J).$
Here ${\underline {\cal B}}$ is the space of unparameterized map 
corresponding to ${\cal B},$ and we obtain ${\underline {\cal B}}$ from
${\cal B}$ by quotienting out the ${\bf R}^1\times S^1$ action of the
domains of its elements. 

\end{lemma}

We non come to the non-transversality of the trivial connecting map.
Note that to define ${\cal B}(\tilde u)$, we fixed a $N'$ and for any
$d'\in D_\epsilon(d)\in {\bf R}^4$, we then define the base connecting
map $\udn$. However in the case that $\{x_-\}=\{x_+\}$, for any ${\tilde
u}\in{\tilde{\cal M}}(x_-, x_+; J)$, the two asymptotic ends $x_-'$ and
$x_+'$ of ${\tilde u}$ are the same as parameterized closed orbits. This
suggests that to define $\udn$ we only need to consider those $d'\in
D_\epsilon(d)$ such that $\theta$-component $d_{+,2}'=d_{-,2}'.$
Therefore, we get a subspace ${\cal B}^0\hookrightarrow{\cal B}$ such
that locally, 
$${\cal B}^0\equiv L^p_{1,\delta}({\tilde u}_d^*T(\tilde M))\times
D_\epsilon^0(d),$$ where
$D_\epsilon^0(d)\hookrightarrow D_\epsilon(d)$ with $d'_{+, 2}=d'_{-, 2}$
 and a Banach bundle ${\cal L}\rightarrow {\cal B}^0 $. Here the fiber $
{\cal L}_{\tilde{v}},$  ${\tilde{v}}\in
{\cal B}^0$ defined as before, and $\bp :{\cal B}^0\rightarrow{\cal L}$
is a section 
such that ${\tilde {\cal M}}(x_-, x_+; J)\cap {\cal B}^0$ is the zero
locus of $\bp$. However the linearized $\bp$-operator at $\tilde u,$
$DF(0)=L=L_1\oplus L_2,$ has index 1 rather than 2. In fact in this case
$L_2$ is still an isomorphism and hence $Ind(L_2)=0,$ but $Ind(L_1)=1.$
Since any element $u$ in ${\tilde {\cal M}}(x_-, x_+; J)$ still has the
$2$-dimensional symmetries, this implies that $\bp :{\cal B}^0\rightarrow
{\cal L}$ is not transversal along ${\tilde{\cal M}}(x_-, x_+; J).$

\section{Additivity of Fredholm Index}

We now come to the question on the additivity of the index under gluing
map. Combining with
the compactness theorem in [L1] with this formula will give the desired
virtual codimension of the boundary components of the moduli space of
connecting pseudo-holomorphic maps in the symplectification. In [L2], we
will prove that the codimension can be realized in the corresponding
virtual moduli cycles. 
 Since we have already written down the precise formula for the
index of $L_1$ , we only need to find how $L_2$ changes under gluing.

Let ${\tilde{\cal M}}(x_-, x, x_+; J)$ be the subspace of 
${\tilde{\cal M}}(x_-, x; J)\times {\tilde {\cal M}}(x, x_+; J)$, which
consists of all pairs $(\tilde u, \tilde w)$ such that
$${\tilde u}:  S^2\setminus\{z_-, z_+'\}  \rightarrow {\tilde{\cal M}},
\quad\quad
{\tilde w}:  S^2\setminus\{z_-', z_+\}\rightarrow{\tilde{\cal M}} 
\quad\quad\mbox{ and }$$
 $$\lim_{s\rightarrow z_+'}{\tilde u}(s,t) =  \lim_{s\rightarrow z_-'}
{\tilde w}(s,t)=x'(t), x'\in \{x\}. $$
Given $(\tilde u, \tilde v)\in {\tilde{\cal M}}(x_-, x, x_+; J)$ and a
gluing parameter $\tau\in{\bf R}^+$, we will glue ${\tilde u}_{N'}$
 and
${\tilde w}_{N'}$ together to get a ${\tilde v}_{N', \tau}={\tilde u}_{N'}
\#_\tau{\tilde w}_{N'}$. Here we use  ${\tilde u}_{N'}$ to denote  ${\tilde
u}_{d,N'}$ defined before with $d\in{\bf R}^{2(l+m)}$ describing the asymptotic limits of
${\tilde u}$. Similarly for  ${\tilde w}_{N'}. $ Note that since there is
no closed loop in the bubble tree, we have assumed in above case that
$z_+'$ and $z_-'$ has only one element. 

Since both ${\tilde u}_{N'}$ and ${\tilde w}_{N'}$ are trivial maps when
$|s|$ is large enough and they have the same asymptotic limits along the
end $z_+'=z_-'$ as parameterized closed orbit, therefore, given
$\tau\in{\bf R}^+$, there is an obvious way to define ${\tilde v}_{N',
\tau}={\tilde u}_{N'}\#_\tau{\tilde w}_{N'}.$ Namely, we cut off the
part of the cylinder $S^2\setminus\{z_-, z_+'\}$ near $z_+' $ with $s>\tau$ and the
corresponding part of the cylinder of $S^2\setminus\{z_-', z_+\}$ hear
$z_-'$ with $s<-\tau$ and glue the remaining parts together along their
boundaries with respect to the same $\theta$-coordinate (the $S^1$-
parametrization along $z_-'$ and $z_+'$, see [L1].) 
To glue the targets
of ${\tilde u}$ and $\tilde w$, we note that when $|s|> N$ near $z_-'$
and $z_+'$, the $a$-parts of $\tilde u$ and $\tilde w$ are linear maps
only depending on $s$. 
Let ${\tilde M}_u$ and ${\tilde M}_w$ be the
targets of $\tilde u$ and $\tilde w$. We get ${\tilde M}={\tilde
M}_v={\tilde M}_u\#_\tau{\tilde M}_w$ by cutting off the parts of
${\tilde M}_u$ with $a >a(\tau)$ and the part of ${\tilde M}_w$ with
$a<a(-\tau)$ and glue the remaining parts together. 

We then define ${\tilde v}_{N', \tau}: S^2_\tau\setminus 
\{z_-, z_+\}\rightarrow {\tilde M}_v$ by simply applying ${\tilde
u}_{N'}$ and ${\tilde w}_{N'}$ to the points in $S^2_\tau\setminus
\{z_-, z_-\}$ which are in the domains of ${\tilde u}$ and $\tilde w.$
Here $S^2_\tau\setminus\{z_-, z_+\}$ is the domain of ${\tilde v}_{N',
\tau}$. It is easy to see that when $\tau$ is large enough, ${\tilde
v}_{N', \tau}$ is well-defined. To simplify our notation, we denote 
${\tilde v}_{N', \tau}$  by ${\tilde v}_\tau.$ 

Now ${\tilde v}_\tau={\tilde v}_{N',\tau}$ plays the same rule as
${\tilde u}_{N'}$ and ${\tilde w}_{N'}$ in ${\tilde{\cal M}}(x_-, x;
J)$ and ${\tilde {\cal M}}(x, x_+; J).$ It is a connecting map between
$x_-$ and $x_+$, which is almost ${\tilde J}$-holomorphic and is the
trivial map near the two ends $x_-$ and $x_+.$ We have three linearized
$\bp$-operator $L_{{\tilde u}_{N'}},$  $L_{{\tilde w}_{N'}}$ and
$L_{{\tilde v}_\tau}.$ Let $D_u$, $D_w$ and $D_{v_\tau}$ be the $L_2$-
part of these operators. Then
$$D_{u,w}=(D_u, D_v):  L^p_{1,\delta}({\tilde u}_{N'}^*\xi)\oplus
L^p_{1,\delta}({\tilde w}_{N'}^*\xi) \rightarrow
L^p_{0,\delta}({\tilde u}_{N'}^*\wedge^{0,1}(\xi))\oplus
L^p_{0,\delta}({\tilde w}_{N'}^*\wedge^{0,1}(\xi)),
$$
and 
$$
D_{v_\tau}:L^p_{1, \delta}({\tilde v}_\tau^*\xi)\rightarrow
L^p_{0,\delta}({\tilde v}_\tau ^*\wedge^{0,1}(\xi)).$$

\begin{pro}

When $\tau$ is large enough, 
$$Ind(D_{v_\tau})=Ind(D_{u,w}).$$
\end{pro}

\proof

As we mentioned before, for the calculation of the indices here, we can
use $L^p$-norm rather than the weighted $L^p$-norm. 

In the case that both domains of $\tilde u$ and $\tilde w$ are the
cylinder (hence the domain of ${\tilde v}_\tau$ is the cylinder too),
the proposition follows from the fact that in this case the index
involved can be calculated by using the spectral flow. The detail of the
argument of this type is given in [F] by Floer. 

For the general case, the proposition can be proved by using a linear
version of the so called 'gluing' technique in Floer homology and
quantum homology. We will outline this argument here. For details  of
the non-linear version of this argument, see [F], [L] and [LT]. 

We assume that $D_{u,w}$ is surjective. Then general case can be
reduced to this case (see [LT]). We will show that when $\tau$ is large
enough, $D_{v_\tau}$ is also surjective and $dim(ker D_{v_\tau})=dim
D_{u,w}.$ For this purpose, we assume that $K_u=ker(D_u)={\bf
R}\{f_1,\cdots, f_m\}$ and $K_w=ker(D_w)={\bf R}\{g_1, \cdots, g_n\}$.
Then $K_{u,w}=ker(D_{u,w})={\bf R}\{f_1, \cdots, f_m, g_1, \cdots,
g_n\}.$ Let $\beta_{\tau, u}, \beta_{\tau, w}$ be the cut-off functions
specified by 
$$
\beta_{\tau, u}(s)=\left \{ \begin{array}{lr}
			 1 & s<\tau-1\\	
			0 & s>\tau
			\end{array}\right. 
\quad \mbox{ and }
\beta_{\tau, w}(s)=\left \{ \begin{array}{lr}
			1 & s>-\tau +1\\
			0 & s< -\tau
			\end{array}\right. .$$	
				
Define $f_{\tau, i}=\beta_{\tau, u}\cdot f_i$ and $g_{\tau,
j}=\beta_{\tau, w}\cdot g_j$. Let 
$$N_\tau={\bf R}\{f_{\tau, 1}, \cdots,
f_{\tau, m}, g_{\tau, 1}, \cdots, g_{\tau, n}\}$$ 
be the vector space generated by the asymptotic kernel of $D_{v_\tau}$.
Note that each $f_{\tau, i}$ and $g_{\tau, j}$ is defined over ${\tilde
v}_\tau. $ Let $K_{\tau}=ker(D_{v_\tau}).$ Let $N_\tau^\perp$ and
$K_\tau^\perp$ be the $L^2$-complement of $N_\tau$ and $K_\tau$ in
$L^p_1({\tilde v}_\tau^*\xi) $ respectively. 

\begin{sublemma}
There exists a constant $C>0$ independent of $\tau$ such that for any
$\eta\in N^\perp_\tau$ with $\tau$ large enough, 
$$\|\eta\|_{1,p}<\|D_{v_\tau}(\eta)\|_{0,p}.\label{1}$$
\end{sublemma}

\proof
Similar statement has been proved in [F], [L] and [LT].  

Assume \ref{1} is not true. Then there exists $\eta_\tau\in
N^\perp_\tau$ such that 

\noindent (a) $\|\eta_\tau\|_{1,p}=1;$ and 

\noindent (b) $\|D_{v_\tau}(\eta_\tau)\|_{0,p}\rightarrow 0,$

\noindent as $\tau\rightarrow\infty.$ 

We want to show that (a) and (b) contradicts to each other. 

Let $(\tilde s, t)$ be the new cylindrical coordinate for the "neck"
part of the domain $S^2_\tau\setminus\{z_-, z_+\}$ of ${\tilde v}_\tau$,
which starts at the middle of the neck. Then ${\tilde s}\in (-\rho_\tau,
\rho_\tau)$, where $\rho_\tau=\tau-N'.$ We want to show that in terms of
this coordinate, $\|\eta_\tau|_{[-2, 2]\times S^1}\|_{1,p}\rightarrow 0$
as $\tau\rightarrow\infty.$ To this end, let ${\tilde \eta}_\tau$ be
$\eta_\tau|_{(-\rho_\tau, \rho_\tau)\times S^1}$ with respect to
$({\tilde s}, t)$-coordinate. Since $\|{\tilde \eta}_\tau\|_{1,p}< 1$,
${\tilde \eta}_\tau$ is weakly convergent to ${\tilde \eta}_\infty\in
L^p_1({\bf R}\times S^1, {\bf R}^{2n})$ in $L^p_1$-norm.
Note that here we have used the fact that in the neck part $(-\rho_\tau, 
\rho_\tau)\times S^1, $ ${\tilde \eta}_\tau\in L^p_1((-\rho_\tau,
\rho_\tau)\times S^1, {\bf R}^{2n})$. One can show that (b) implies that
$D_\infty{\tilde \eta}_\infty =0. $ Here $D_\infty =\cup_\tau
D_{v_\tau}$. Note that if $\tau <\tau ', D_{v_{ \tau '}}|_{(-\rho_\tau,
\rho_\tau)\times S^1}=D_{v_\tau}.$
However $D_\infty$ is an isomorphism. We have ${\tilde \eta}_\infty =0$.
It follows from Sobelev embedding theorem that on any compact subset $[-
R, R]\times S^1$, 
${\tilde \eta}_\tau|_{[-R, R]\times S^1}\rightarrow 0$ in $C^0$-norm. Since 
$\|D_{v_\tau}{\tilde \eta}_\tau\|_{0,p}\rightarrow 0$, standard elliptic estimate
implies that $\|{\tilde \eta}_\tau|_{[-2,2]\times S^1}\|_{1,p}\rightarrow 0$ as
$\tau\rightarrow \infty.$ 

Let $$\beta(\tilde s)=\left \{\begin{array}{lr}
			1 & -2 <{\tilde s}< 2\\
			0 & \mbox{ otherwise }
				\end{array}\right.
$$
be a cut-off function.  Write $\eta_\tau=(1-\beta)\eta_\tau
+\beta\eta_\tau.$ Then $D_{v_\tau}((1-\beta)\eta_\tau)=D_{u,w}((1-
\beta)\eta).$ Hence 
\begin{eqnarray*}
\|(1-\beta)\eta_\tau\|_{1,p} & \leq  & C\|D_{u,w}(1-\beta)\eta_\tau\|_{0,p}
+|\pi (1-\beta)\eta_\tau|\\
& = & C\|D_{v_\tau}((1-\beta)\eta_{\tau})\|_{0,p}+|\pi(1-
\beta)\eta_\tau|\rightarrow 0.
\end{eqnarray*}
Here $\pi$ is the orthogonal projection map to the kernel of $D_{u,w}.$
Therefore, $\|\eta_\tau\|_{1,p}\leq \|(1-
\beta)\eta_\tau\|_{1,p}+\|\beta\eta_\tau\|_{1,p}\rightarrow 0$ as
$\tau\rightarrow \infty.$ 

Therefore(a) and (b) contradict to each other.
\QED

Now define $${\tilde D}_\tau:L^p_1({\tilde v}^*_\tau\xi)\rightarrow
L^p_0({\tilde v}_\tau ^*\wedge^{0,1}(\xi))\oplus N_\tau$$
 and
$${\tilde D}_{u,w}: L^p_1({\tilde u}^*_{N'}\xi)\oplus L^p_1({\tilde
w}^*_{N'}\xi)\rightarrow L^p_0({\tilde u}^*_{N'}\wedge^{0,1}(\xi))\oplus
L^p_0({\tilde w}^*_{N'}\wedge^{0,1}(\xi))\oplus K_{u,w}$$ as follows. 

Define ${\tilde D}_{u,w}=D_{u,w}\oplus \pi_K$, where $\pi_K$ is the
orthogonal projection. Similarly, ${\tilde
D}_\tau=D_{v_\tau}\oplus\pi_\tau$, where $\pi_\tau$ is the projection to
$N_\tau .$ 

The previous sublemma shows that ${\tilde
D}_\tau|_{N^\perp_\tau}=D_{v_\tau}|_{N^\perp_\tau}$ is injective.
Therefore, ${\tilde D}_\tau$ is injective. Note that $dim N_\tau=dim
K_{u,w}$. Therefore we only need to show that ${\tilde D}_\tau$ is also
surjective when $\tau$ is large enough (note that ${\tilde D}_{u,w}$ is
surjective, which follows from our assumption). To this end, 
it is more convenient to use $L^2_2$-norm for the domain of $D_{u,w}$
and $D_{\tau}$ before. 
Observe that
since each element $\phi\in ker_{u,w}$ or $N_\tau$ is in $L^2_2$, the
function $<\phi, >$ is continuous on $L^2_{-2}$.
  Therefore, $\pi_K$ and  $\pi_\tau$ and
hence ${\tilde D}_{u,w}$ and ${\tilde D}_\tau$ are well-defined by
replacing $L^2_2$-norm of the domain with $L^2_{-2}$-norm, and $L^2_1$-
norm od the target  with $L^2_{-3}$-norm. Now passing to the dual of ${\tilde
D}^*_\tau$ and ${\tilde D}^*_{u,w}$ and note that the domains of the
operators have $L^2_3$-norm and the ranges have $L^2_2$-norm. 

We now in the position to apply the previous lemma to deduce that the
injectivity of ${\tilde D}^*_{u,w} $ implies the injectivity of ${\tilde
D}^*_\tau$ when $\tau $ is large enough.
\QED

Now we can compare the virtual dimension of ${\tilde{\cal M}}(x_-, x,
x_+; J)$ near $({\tilde u}_{N'}, {\tilde w}_{N'})$ with the dimension of
${\tilde {\cal M}}(x_-, x_+; J)$ near ${\tilde v}_\tau.$ We start with
 two particular but most important cases. The first case is that
${\tilde u}_{N'}\in {\tilde{\cal M}}(x, J)$ is an almost holomorphic
plane with a positive end $x$ and ${\tilde w}_{N'}\in{\tilde {\cal
M}}(x_-, x, x_+; J)$ is a connecting almost holomorphic map of three
ends $(x_-, x)$ and $x_+$, where $(x_-, x)$ are negative ends and $x_+$
is the positive one such that their asymptotic limits as parameterized
closed orbits agree at $x$. Let ${\tilde{\cal M}}(x; (x_-, x), x_+; J)$
be the moduli space of all such pairs. This corresponding to the case of
only one bubble. As we mentioned before, in the simplest case of
of  bubbling off only one bubble, the domain splits into three 
components rather than just two as considered here. However, as far as 
the computation of the index is concerned, we still can consider above
case first.
 
We want to compare the local dimension,
$dim_{{loc}_{(u,w)}}\{{\tilde{\cal M}}(x; (x_-, x), x_+; J)\}$ near
$({\tilde u}_{N'}, {\tilde w}_{N'})$ with the local dimension
$dim_{{loc}_{v_\tau}}{\tilde{\cal M}}(x_-, x_+; J)$ near
${\tilde v}_\tau.$
Note that ${\tilde{\cal M}}(x_-, x_+)$ is the moduli space of connecting
maps between these two ends $x_-$ and $x_+$. 

Now 
\begin{eqnarray*}
& & dim_{{loc}_{(u,w)}}\{{\tilde{\cal M}}(x; (x_-, x), x_+; J)\} \\
& & = (Ind(L_{1, u})+ Ind(L_{1, w})-1)+Ind(D_{u,w})\\
& & =  (2+2-1)+Ind(D_{u,w})=2+Ind(D_{v_\tau})+1 \\
& &  =  Ind(L_{1, {v_\tau}})+Ind(D_{v_\tau})+1 \\
& & = dim_{{loc}_{v_\tau}}\{{\tilde{\cal M}}(x_-, x_+; J)\}+1. 
\end{eqnarray*}

Now the dimension of the symmetries of the domain of $({\tilde u}_{N'},
{\tilde w}_{N'})$ is 3 while the dimension of the symmetries of the
domain of ${\tilde v}_\tau$ is 2. Hence for unparameterized curves, the
two moduli spaces have the same dimension. Finally, put the ${\bf R}$-
symmetry in target into consideration, we have:
\begin{pro}

$dim_{{loc}, (u, w)}({\cal M}(x; (x_-, x), x_+; J))=dim_{loc,
v_\tau}({\cal M}(x_-, x_+; J))-1.$

That is bubbling off only one bubble is a codimension one phenomenon.
\end{pro}

The second important case is that both ${\tilde u}_{N'}$ and ${\tilde
w}_{N'}$ are connecting maps between two ends. As above, a direct
calculation shows the following:

\begin{pro}

$dim_{loc, (u, w)}({\cal M}(x_-, x, x_+; J))=dim_{loc, v_\tau}({\cal
M}(x_-, x_+; J))-1.$
That is splitting a connecting map into broken ones of two elements is
also a codimension one phenomenon.

\end{pro}

A direct calculation shows that both of above statements are true also
for the multi-ends case. We have:
\begin{pro}

For a family of connecting pseudo-holomorphic maps with multi-ends, both
bubbling off one bubble and splitting off a connecting map of two ends
are codimension one  phenomenon.
 
\end{pro}

Note that in the case that the above connecting maps  
split off a trivial connecting  map, we have used the reduced Banach space
neighborhood ${\cal B}^0(u)$ for the trivial map to compute the relevant
 index here.

\begin{pro}

When a family of connecting pseudo-holomorphic map of multi-ends splits
into a family of broken connecting pseudo-holomorphic maps of two
elements, each being multi-ends, the splitting has codimension two. 
\end{pro}
\proof
Here "multi-ends" means that the number of the ends here is at least
three.
The proof follows from a direct computation of the relevant indices. 
\QED
Not that in the last two statements, we allow the marked points to move. 

Finally, combine these proposition together,
we have the following theorem

\begin{theorem}

Given a family of connecting pseudo-holomorphic maps, virtually, 
the bubbling
is a co-dimension one phenomenon while the splitting off the connecting
maps into the broken ones of two elements is codimension one.

\end{theorem}

\end{document}